\documentclass[12pt]{article}
\usepackage{amssymb,amsmath}

\newtheorem{thm}     {Theorem}[section]
\newtheorem{prop}    [thm]{Proposition}

\newtheorem{cor}     [thm]{Corollary}
\newtheorem{lemma}   [thm]{Lemma}
\newtheorem{remark}   [thm]{Remark}

\newcommand{\proof} {\noindent{\bf Proof. }}
\newcommand{\E}{{\mathcal E}}
\newcommand{\N}{\mathbb N}
\newcommand{\R}{\mathbb R}
\def\Re{{\rm Re\,}}
\def\TBP{{\rm TBP} }
\def\FP{{\rm FP}}
\def\T{{\rm T}}
\def\rank{{\rm rank}}
\def\tr{{\rm tr}}
\def\t{\theta}

\begin{document}

\title{Minimal biquadratic energy of 5 particles on 2-sphere}
\author{Alexander Tumanov}
\date{}
\maketitle

{\small
University of Illinois, Department of Mathematics,
1409 West Green Street, Urbana, IL 61801, USA,
tumanov@illinois.edu
}
\bigskip
\bigskip
\bigskip

Abstract.
Consider $n$ points on the unit 2-sphere. The potential energy
of the interaction of two points is a function $f(r)$ of the
distance $r$ between the points. The total energy $\E$ of $n$
points is the sum of the pairwise energies. The question is how
to place the points on the sphere to minimize the energy $\E$.
For the Coulomb potential $f(r)=1/r$, the problem goes back to
Thomson (1904). The results for $n < 5$ are simple and well known.
We focus on the case $n=5$, which turns out to be difficult.
In this case, the following results have been obtained.
Dragnev, Legg, and Townsend \cite{DLT} give a solution
of the problem for $f(r)= -\log r$ known as Whyte's problem.
Hou and Shao give a rigorous computer-aided solution for $f(r)= -r$.
Schwartz \cite{Sch} gives a rigorous computer-aided solution
of Thomson's problem. We give a solution for biquadratic
potentials.
\bigskip


\bigskip

\section{Introduction}

Let $P=\{p_1,...,p_n\}$ be a configuration of $n$ points
on the unit sphere $S^2\subset\R^3$.
Let $f$ be a real function that we regard as the potential
of the interaction of two points.
Then the potential energy of $P$ has the form
\[
\E_f(P)=\sum_{i<j}f(|p_i-p_j|),
\]
here $|x|$ denotes the standard Euclidean norm of $x\in\R^3$.
The problem of minimizing $\E_f(P)$ is a subject of
hundreds of publications
(see, e. g., \cite{CK} and references there).
For the Coulomb potential $f(r)=1/r$, the problem has
a significance in physics and goes back to Thomson (1904).

We are concerned with the energy minimizing problem
for small values of $n$.
To eliminate trivial or irrelevant potentials,
we assume that $f$ is decreasing and strongly convex
as a function of $r^2$. That is, for the function
$\phi(t)= f(\sqrt{t})$, we have $\phi'<0$ and $\phi''>0$.
The results for $n < 5$ are simple and well-known.
They are special cases of a more general
result that the regular simplex is a minimizing
configuration in $S^{m-1}\subset\R^m$ for every dimension $m$
(see, e. g., \cite{CK}).

We focus on the case $n=5$, which turns out to be difficult.
In this case, the following results have been obtained.
Dragnev, Legg, and Townsend \cite{DLT} give a solution of the
problem for $f(r)= -\log r$ known as Whyte's problem.
Hou and Shao \cite{HouShao} give a rigorous computer-aided
solution for $f(r)= -r$, for which the problem is well-known
in discrete geometry. Schwartz \cite{Sch} gives a rigorous
computer-aided solution of Thomson's problem.
The results of \cite{HouShao} and \cite{Sch}
involve massive calculations that require a computer.
In all these three results, a unique minimizer is the so-called
triangular bipyramid (TBP) that consists of two antipodal points,
say the North and South poles, and three points on the equator
forming an equilateral triangle.
However, following \cite{CK} and \cite{Sch} we point out that
TBP is not universally optimal for either all inverse power
potentials $f(r)=r^{-a}, a>0$, or potentials of the form
$f(r)=(4-r^2)^k, k\in\N$. For sufficiently large values of
$a$ and $k$, the energy takes a smaller value on another
configuration, a regular four-pyramid with square base (FP).
The latter depends on a parameter, say the pyramid's height,
whose value is chosen to minimize the energy and it depends on $f$.

Clearly, a minimizing configuration must be a critical point
for the energy $\E_f$. We conjecture that for $n=5$ and
almost all potentials $f$, there are only two nontrivial critical
configurations: TBP and FP.
Here we call a configuration trivial if it has repeated points
or reduces to lower dimension.
It would suffice to prove the conjecture for potentials $f$
which are polynomials in $r^2$.
The case of $f(r)=-r^2$ is degenerate: it is not strongly
convex in $r^2$.
In this case the answer is the following: a configuration $P$
is a minimizer exactly if $\sum p_i=0$, see Lemma \ref{pre1} below.
We consider the first non-degenerate case $f(r)=ar^4-br^2+c$,
which we call a {\it biquadratic} potential.
This function is decreasing and strongly convex in $r^2$ if
$a>0$ and $b>8a$.
This potential does not have a clear geometric or physical
meaning. However, from the point of numerical calculations,
it is similar to the physical potential $1/r$.
Indeed, simple analysis (see \cite{Sch}, Lemma 2.1)
shows that the minimizer for $1/r$ cannot have points
within the distance $1/2$.
On the rest of the interval, that is, for
$1/2\le r\le 2$, for suitable coefficients $a, b$, and $c$,
the biquadratic polynomial approximates $1/r$ to high accuracy.
Our main result is the following.

\begin{thm}
\label{MainThm}
Let n=5 and let $f(r)=ar^4-br^2+c$, $a>0$, $b>8a$.
Then \TBP and \FP are the only two nontrivial critical
configurations for $\E_f$ and \TBP is a unique global
minimizer.
\end{thm}

One motivation for this result is of course to obtain
a solution for a non-trivial potential without massive
computer-aided calculations. We point out that known
solutions for several potentials may yield a solution
for another potential. Indeed, suppose that for several
potentials $f_i$, it is known that TBP is a minimizer.
Let $g=\sum c_i f_i$, $c_i>0$, approximate a potential
$f$ in such a way that $g(r)\le f(r)$ for all $0<r\le 2$,
and $g$ matches with $f$ on all distances $r$ between the
points of TBP, that is, $r=\sqrt{2}, \sqrt{3}$, and $2$.
Then TBP is a minimizer for $f$. For instance, if TBP is
known to be a minimizer for $(4-r^2)^k$ for $k=0,1,2,3,5$,
then using the above argument one can show that TBP will
be a minimizer for both power potentials
$r^{-a}$ and $-r^a$ for all $0<a\le 2$. Our biquadratic
potential covers $k=0,1,2$, but we admit that proving the
result for $k=3$ and 5 would be hard.

We now describe our strategy for proving the main result.
We start Section 2 by rewriting the energy in terms
of the matrix of inner products $b_{ij}=p_i\cdot p_j$.
Then the original problem reduces to minimizing
$E_h(P)=\sum_{i<j} h(b_{ij})$, where $h(t)=(t+a)^2$.
We then show that it suffices to consider only
one function $h(t)=(t+1)^2$. We also prove in this section
that a minimizer cannot have repeated points.

In Section 3 we introduce equations for critical
points of the energy, that is, the equilibrium equations.
We use a special orthogonal coordinate system such that
the matrix whose rows represent the points $p_i$
has orthogonal columns. This coordinate system
significantly simplifies both the biquadratic energy
and equilibrium equations.

In Section 4, using complex analysis language, we obtain
an auxiliary result on non-degeneracy
of a Cauchy type matrix (Theorem \ref{Cauchy}) that may
be of interest by itself. We use it to eliminate the
bulk of non-critical configurations.

In Section 5 we use the obtained simple form of
the equilibrium equations and a case by case analysis
to reduce the problem to what we call Main Special Case.
The latter restricts the problem to configurations having
three points on a great circle and two other points symmetric
about that circle.

In Section 6 we directly solve the system in Main Special Case.
Surprisingly, even in this simple case
the system gives us enough trouble that in the end we use
symbolic computation software. We point out however that
the use of a computer is light and can be completely eliminated.

\section{Preliminaries}

For the sake of generality and simplicity of notations,
we sometimes consider $n$ points on the unit sphere
$S^{m-1}\subset \R^m$ without specifying the values of $m$ and $n$.
Our main result requires $m=3$ and $n=5$, however, some other
results of the paper hold for general $m$ and $n$.
The reader may always assume without much loss that
$m=3$ and $n=5$.

We first rewrite the energy in terms of the pairwise
inner products. Let $x,y\in S^{m-1}$, $m\in\N$.
Let $t=x\cdot y$ be the Euclidean inner product,
and let $r=|x-y|$. Then $r^2=2-2t$.
For a function $f(r)$, we define $h(t)=f(\sqrt{2-2t})$.
If $f$ is decreasing and convex on $(0,2]$, then $h$ is
increasing and convex on $[-1,1)$. Then we rewrite the energy
$\E_f(P)$ of a configuration
$P=\{p_1,\dots,p_n\}$ in $S^{m-1}$ in the form
\begin{equation}
\label{Eh}
E_h(P)=\sum_{i<j} h(p_i\cdot p_j).
\end{equation}
For a biquadratic potential after omitting unnecessary
constants we have
\[
h(t)=(t+a)^2, \quad
a\ge 1.
\]
Note that for $a\ge1$ the function $h$ is increasing.
We recall the following simple fact.
\begin{lemma}
\label{pre1}
Let $h(t)=t$. Then $P=\{p_1,\dots,p_n\}$ is a minimizer
for $E_h$ if and only if $\sum p_i=0$.
\end{lemma}

\proof
Follows from
$E_h(P)=\sum_{i< j} p_i\cdot p_j={1\over 2}(|\sum_{i=1}^n p_i|^2-n)$.
$\square$

\begin{cor}
\label{pre2}
Let $n=5$, $m=3$.
It suffices to show that \TBP is a unique minimizer for
$E_h$, $h(t)=(t+1)^2$. Then it would follow that
\TBP is a unique minimizer for $E_h$, $h(t)=(t+a)^2$,
for all $a\ge 1$.
\end{cor}

\proof
Note that $h(t)=(t+a)^2$ can be written as
$h(t)=h_0(t)+h_1(t)$, here $h_0(t)=(t+1)^2$,
$h_1(t)=C_1 t+C_2$, $C_1=2(a-1)\ge0$, $C_2=a^2-1$.
Since for \TBP we have $\sum p_i=0$, then by Lemma \ref{pre1},
\TBP is a minimizer for $E_{h_1}$.
Now since \TBP is a unique minimizer for $E_{h_0}$, then
it is also a unique minimizer for $E_h$.
$\square$
\medskip

For inverse power potentials $f(r)=r^{-\alpha}$,
a minimizer can not include repeated points because
$f(0)=\infty$. The latter is not the case for biquadratic
potentials. For convenience we consider the case of
repeated points separately.

\begin{prop}
\label{pre3}
A minimizer for the biquadratic energy
$E_h$, $h(t)=(t+1)^2$, of 5 points on $S^2$ can not have
repeated points.
\end{prop}

\proof
We use the notation $E(P)=E_h(P)$ for $h(t)=(t+1)^2$.
For TBP the inner products are the following:
one $-1$, six 0, and three $-1/2$.
Hence $E(\TBP)=0+6+3(1/2)^2=27/4$.

We also use a regular tetrahedron T on the sphere.
We recall that T is a universal minimizer \cite{CK},
in particular, for our potential.
For T, all the four inner products equal $-1/3$.
Hence $E(\T)=6 (2/3)^2=8/3$.

Arguing by contradiction, suppose
$P=\{p,p,q_1,q_2,q_3\}$ is a minimizer.
Then we have
\[
E(\TBP)\ge
E(P)=E(\{p,p\})+E(\{ p,q_1,q_2,q_3 \})+
\sum_{i=1}^3 E(\{p,q_i\}).
\]
Since
$E(\{ p,q_1,q_2,q_3 \})\ge E(\T)$, then
\[
\sum_{i=1}^3 E(\{p,q_i\})\le
E(\TBP)-E(\T)-E(\{p,p\})=27/4-8/3-4=1/12.
\]
Hence $E(\{p,q_i\})\le 1/12$, $i=1,2,3$.
Put $t=p\cdot q_i$. Then $(t+1)^2\le 1/12$,
and $t\le \sqrt{3}/6-1 < -\sqrt{2}/2$.
Since the angle between $-p$ and $q_i$ is less
that $\pi/4$, then the angle between
$q_i$ and $q_j$ is less than $\pi/2$, that is,
$q_i\cdot q_j>0$, and $E(\{q_i,q_j\})>1$.
We now have
\[
27/4=E(\TBP)\ge E(P)\ge
E(\{p,p\})+E(\{q_1,q_2,q_3 \})\ge
4+3=7,
\]
which is absurd.
$\square$

\section{Equilibrium equations}

We find the minimizer by finding all critical
points of the energy functional. We call them critical
or equilibrium configurations.

Let $P=\{p_1,...,p_n\}\subset S^{m-1}$.
Put $b_{ij}=p_i\cdot p_j$.
We apply the Lagrange multipliers method to minimizing
the energy (\ref{Eh}) subject to $n$ constraints $|p_i|=1$.
By differentiating the Lagrange function
\[
L(P,\mu)=E_h(P)-\sum_{i=1}^n \mu_i\, |p_i|^2
\]
with respect to all components of $p_i$
we obtain the equation
\[
\sum_{j\ne i}h'(b_{ij})p_j=2\mu_i p_i.
\]
This equation physically means that the resulting
force acting on $p_i$ from the rest of the particles
is orthogonal to the sphere.
By inner-multiplying both sides by $p_i$, we can solve for
$\mu_i$. After eliminating $\mu_i$ we obtain
a system of $n$ vector equations for critical configurations
\begin{equation}
\label{Crit1}
\sum_{j=1}^n h'(b_{ij})(p_j-b_{ij} p_i)=0,\quad
1\le i\le n.
\end{equation}
Note that if $h'$ is continuous at 1,
which is the case for biquadratic potentials,
then we do not have to
omit the term with $j=i$ because it is equal to 0.
The system is difficult to solve in general.
Following \cite{DLT}, we point out that
in the case of the logarithmic potential
$h(t)=-\log(1-t)$, by adding up all the equations,
one obtains $\sum p_i=0$, which is useful
in solving the system (\ref{Crit1}).
However, for the biquadratic potential, this method does not
work. Although the equation $\sum p_i=0$ is valid
for \TBP, it is not valid for all critical configurations,
in particular for \FP.

Note that the system is invariant under rotations
of the sphere. To work with independent parameters,
one needs a normalization. A natural normalization for $m=3$
would consist of fixing one point, say at the North pole,
and restricting position of the second point to a fixed
meridian. However, we use a different normalization.

Let $X$ be the matrix whose rows are the points
$p_i$ of the configuration $P$. Our normalization
requires that the matrix $X$ have orthogonal columns.
We show that it is always the case in a suitable
orthogonal coordinate system in $\R^m$.

For a matrix $A$, let $A^*$ denote the transpose of $A$.

\begin{lemma}
\label{equilib1}
Let $X$ be a real $n\times m$ matrix.
Then there is an orthogonal $m\times m$ matrix $U$
such that the matrix $Y=XU$ has orthogonal columns,
that is, $Y^*Y$ is diagonal.
\end{lemma}

\proof
We have $Y^*Y=(XU)^*(XU)=U^*AU$, here $A=X^*X$.
The matrix $A$ is a symmetric $m\times m$ matrix,
hence it is diagonalizable in a suitable orthonormal basis.
That is, there is an orthogonal matrix $U$ such that
$\Lambda=U^*AU$ is diagonal as desired.
$\square$
\medskip

Note that in the above proof, the diagonal of $\Lambda$
consists of the eigenvalues of $A$. Also note
$\rank(\Lambda)=\rank(A)=\rank(X)\le m$.
Finally, we obtain the following.

\begin{cor}
\label{equilib2}
Let $X$ be a real $n\times m$ matrix,
such that $\Lambda=X^*X$ is diagonal.
Let $B=XX^*$ and let $\lambda_1,...,\lambda_m$ be
the diagonal entries of $\Lambda$. Then the columns of
$X$ are eigenvectors of $B$ with the eigenvalues
$\lambda_1,...,\lambda_m$, that is $BX=X\Lambda$.
The other eigenvalues of $B$ are zeros.
\end{cor}

Let $P=\{p_1,...,p_n\}\subset S^{m-1}$ satisfy the above
normalization. Let $X=(x_{ik})$ be the matrix of the coordinates
of all $p_i$. Let $B=(b_{ij})=XX^*$ have the eigenvalues
$\lambda_1,...,\lambda_m$. Then we have
\begin{equation}
\label{Crit2}
\sum_{i=1}^n x_{ik}x_{il}=\lambda_k\delta_{kl}, \quad
\sum_{j=1}^n b_{ij}x_{jk}=\lambda_k x_{ik}.
\end{equation}
Here $k$ and $l$ run from 1 to $m$, $\delta_{kl}$
is the Kronecker symbol.
Since $b_{ii}=|p_i|^2=1$, then $\tr (B)=n$.
Since $\tr (XX^*)=\tr (X^*X)=\sum \lambda_k$, then
\[
\sum_{k=1}^m \lambda_k=n.
\]

We now apply our normalization to the biquadratic energy
and the equilibrium equations.
Denote $\bar x=\sum p_i$.
Taking into account that
$\sum_{ij}b_{ij}^2=\tr(B^2)=\sum\lambda_k^2$, we obtain
\begin{equation}
\label{Crit3}
E(P)=\sum_{i<j}(b_{ij}+1)^2
={1\over 2}\sum_{i,j=1}^n(b_{ij}+1)^2-2n
={1\over 2}\sum_{k=1}^m\lambda_k^2 +|\bar x|^2 +{n^2\over 2}-2n.
\end{equation}
This expression would attain its minimum if $\bar x=0$
and $\lambda_1=...=\lambda_m=n/m$.
However, as we will see later (Proposition \ref{equal-lambdas}),
for $m=3$ and $n=5$ this situation cannot occur, but it can
occur for $m=3$, $n>5$.
This observation makes the case $n=5$ most interesting.

For $m=3$ and $n=5$ we use \eqref{Crit3} to eliminate
configurations in lower dimension.
\begin{prop}
\label{equilib3}
A minimizer for the biquadratic energy
of 5 points on $S^2$ can not lie in a 2-plane.
\end{prop}

\proof
Let $P$ be a critical configuration lying in a plane.
Then the plane must pass through the origin,
otherwise $P$ will not be in equilibrium.
By plugging $n=5$, $m=2$, $\bar x=0$ and $\lambda_1=\lambda_2=5/2$
in \eqref{Crit3} we obtain $E(P)\ge 35/4$, which is greater
than $E(\TBP)=27/4$.
$\square$
\medskip

For the biquadratic potential $h(t)=(t+1)^2$ the system
(\ref{Crit1}) takes the form
\begin{equation}
\label{Crit4}
\sum_{j=1}^n (b_{ij}+1)(x_{jk}-b_{ij}x_{ik})=0.
\end{equation}
Denote
\[
\bar x=(\bar x_1,..., \bar x_m)=\sum_{i=1}^n p_i,\quad
\alpha_i=\sum_{j=1}^n (b_{ij}+b_{ij}^2).
\]
Then by \eqref{Crit2} the system (\ref{Crit4}) further
reduces to
\begin{equation}
\label{Crit5}
(\alpha_i-\lambda_k)x_{ik}=\bar x_k, \quad
1\le i \le n, \quad
1\le k \le m.
\end{equation}
Suppose $\bar x_k\ne 0$ for some $1\le k\le m$.
Then by (\ref{Crit5}) for  all $i$ we have
\begin{equation}
\label{Crit6}
x_{ik}=\frac{\bar x_k}{\alpha_i-\lambda_k}.
\end{equation}
After summing in $i$ and dividing by $\bar x_k$, we have
\begin{equation}
\label{Crit7}
\sum_{i=1}^n \frac{1}{\alpha_i-\lambda_k}=1.
\end{equation}
On the other hand, if $\bar x_k= 0$, then for every $i$
we have either $\alpha_i=\lambda_k$ or $x_{ik}=0$.

As we can see, it matters whether or not $\bar x_k=0$.
Accordingly, we divide the analysis of the system into
cases depending on how many zeros there are among the numbers
$\bar x_k$. We use the notation {\it Case $N$} if there are exactly
$N$ zeros.

\section{Case 0 and Cauchy type matrix}

We recall the Cauchy matrix $C=(c_{ik})$,
$c_{ik}=\frac{1}{\alpha_i-\lambda_k}$ with
distinct $\alpha_i$ and $\lambda_k$.
If $C$ is a square matrix, then it is known
to be nonsingular, that is, $\det C\ne0$.
We need a similar result here.
For future references we include a slightly
more general result than we need.

\begin{thm}
\label{Cauchy}
Let $1\le m\le n$.
Let $\alpha_i$, $1\le i\le n$, and $\lambda_k$,
$1\le k\le m$, be distinct complex numbers satisfying
(\ref{Crit7}). Let $A=(a_{ik})$ be the matrix with entries
\[
a_{ik}=\begin{cases}
\frac{1}{(\alpha_i-\lambda_k)^2} &\text{if}\;\; k\le m,\\
\alpha_i^{k-m-1} &\text{if}\;\; k>m.
\end{cases}
\]
Then $\det A\ne 0$.
\end{thm}
\proof
Suppose there are numbers $r_i$ such that
$\sum_{i=1}^n r_i a_{ik}=0$ for all $1\le k\le n$.
We will prove all $r_i=0$.
Introduce
\[
f(\lambda)=\sum_{i=1}^n \frac{r_i}{(\alpha_i-\lambda)^2},\quad
g(\lambda)=\sum_{i=1}^n \frac{1}{\alpha_i-\lambda}-1.
\]
Then $f(\lambda_k)=g(\lambda_k)=0$ for $1\le k\le m$.
For $k>m$, we have
\begin{equation}
\label{Cauchyeq2}
\sum_{i=1}^n r_i\alpha_i^{k-m-1}=0.
\end{equation}
Note that $g=0$ reduces to an algebraic equation
of degree $n$, hence it has $n$ roots counting multiplicities.
The quotient $h=f/g$ is a rational function vanishing at
infinity. It may have first order poles at $\alpha_i$,
and it may have poles at the zeros $\mu_j$ of $g$ other
than $\lambda_k$.
Hence there exist numbers $s_i$ and $c_{jl}$ such that
\begin{equation}
\label{Cauchyeq3}
h(\lambda)=\frac{f(\lambda)}{g(\lambda)}
=\sum_{i=1}^n \frac{s_i}{\alpha_i-\lambda}
+\sum_j \sum_{l=1}^{\nu_j}\frac{c_{jl}}{(\mu_j-\lambda)^l}.
\end{equation}
Here $\nu_j$ is the multiplicity of $\mu_j$ as a zero of $g$
if $\mu_j$ is not among the numbers $\lambda_k$, and
one unit less otherwise. Then $\sum \nu_j=n-m$.

By passing to the limit as $\lambda\to\alpha_i$ in the
equation \eqref{Cauchyeq3}, we obtain
$s_i=r_i$. We will show that all $c_{jl}=0$.
Denote
$\phi(\lambda)=\sum_{i=1}^n \frac{r_i}{\alpha_i-\lambda}$.
Let $p=n-m$.
The Laurent expansion of $\phi$ at infinity has the form
\[
-\phi(\lambda)=\lambda^{-1}\sum r_i
+\lambda^{-2}\sum r_i\alpha_i+...
+\lambda^{-p}\sum r_i\alpha_i^{p-1}+O(\lambda^{-p-1}).
\]
By \eqref{Cauchyeq2} we obtain
$\phi(\lambda)=O(\lambda^{-p-1})$ as $\lambda\to\infty$.
Note that $\phi'=f$. Then $f(\lambda)=O(\lambda^{-p-2})$,
and in turn $h(\lambda)=O(\lambda^{-p-2})$.
By \eqref{Cauchyeq3} we have
\begin{equation}
\label{Cauchyeq4}
\psi(\lambda)
:=\sum_j \sum_{l=1}^{\nu_j}\frac{c_{jl}}{(\mu_j-\lambda)^l}
=O(\lambda^{-p-1}).
\end{equation}
By reducing to common denominator
$R(\lambda)=\prod (\mu_j-\lambda)^{\nu_j}$,
$\deg R=\sum\nu_j=p$, we obtain
$\psi=Q/R$, where $Q$ is a polynomial.
Then by \eqref{Cauchyeq4},
$Q(\lambda)=O(\lambda^{-1})$. Hence $\psi=0$,
and $h=\phi$.

Since $h'=f$, then the function $h$ satisfies
the differential equation $h'=gh$. Solving this equation
yields
\[
h(\lambda)=Ce^{\int g(\lambda)\,d\lambda}
=Ce^{-\lambda}\prod_{i=1}^n \frac{1}{\lambda-\alpha_i},
\]
which is rational only if $C=0$. Hence $h\equiv 0$ and
all $r_i=0$ as desired. The proof is complete.

In conclusion we note that the proof simplifies a little
for real $\alpha_i$ because in this case all the zeros
of $g$ are simple. The proof further simplifies
in the case $n=m$.
$\square$
\medskip

In the case $m<n$, since the columns of $A$ are linearly
independent, and one of them consists of units,
then we obtain the following.

\begin{cor}
\label{Cauchy2}
In the assumptions of Lemma \ref{Cauchy}, let $m<n$.
Then there do not exist numbers $c_k$, $1\le k\le m$,
such that
$\sum_{k=1}^m \frac{c_k}{(\alpha_i-\lambda_k)^2}=1$
for all $1\le i\le n$.
\end{cor}

We now turn to Case 0.
\begin{prop}
\label{Case0}
Let $m<n$. Then in Case 0 there are no critical
configurations without repeated points.
\end{prop}
\proof
Let $X=(x_{ik})$ be the matrix
of a critical configuration $P$ without repeated points.
Let $\bar x_k\ne 0$ for all $1\le k \le m$.
Then the equation \eqref{Crit6} holds for all $i$ and $k$,
and \eqref{Crit7} holds for all $k$.
If $\alpha_i=\alpha_j$ for some $i\ne j$, then by \eqref{Crit6}
we would have $p_i=p_j$. Hence the numbers $\alpha_i$ are
distinct. If $\lambda_k=\lambda_l$ for some $k\ne l$, then
$\sum_i x_{ik}x_{il}=\sum_i
\frac{\bar x_k \bar x_l}{(\alpha_i-\lambda_k)^2}\ne 0$,
so the columns of $X$ would not be orthogonal.
Hence the numbers $\lambda_k$ are also distinct.
By \eqref{Crit7}, the numbers $\alpha_i$ and $\lambda_k$
satisfy the hypotheses of Corollary \ref{Cauchy2}.
Since $|p_i|=1$, then by \eqref{Crit6} we have
$\sum_k \frac{\bar x_k^2}{(\alpha_i-\lambda_k)^2}=1$.
But then $c_k=\bar x_k^2$ satisfy the equation
in Corollary \ref{Cauchy2}, which is not possible.
$\square$

\begin{remark}{\rm
In the case of our main interest
$m=3$, $n=5$, the first three columns of the matrix
$A$ contain the numbers
$(\alpha_i-\lambda_k)^{-2}$ obtained from
our configuration $P$, the forth column
consists of units, and the fifth column is
unimportant as long as $\det A\ne0$ because all we need is that the
first four columns be linearly independent.
}\end{remark}

\section{Cases 1, 2, and 3}

We restrict to $n=5$, $m=3$.
We look for nontrivial critical configurations.
We regard configurations with repeated points
and those lying in a plane as trivial.
In addition to Cases 0 to 3 we introduce
{\it Main Special Case}, in which three points
lie on a great circle, call it equator, and
the other two are symmetric about the equator.
We will come to the following conclusion.

\begin{prop}
\label{Case123}
All nontrivial critical configurations
in Cases 1, 2, and 3 fall into Main Special Case.
\end{prop}

We proceed with a case by case analysis based on
the equation (\ref{Crit5}), that is,
\[
(\alpha_i-\lambda_k)x_{ik}=\bar x_k, \quad
1\le i \le 5, \quad
1\le k \le 3.
\]

{\bf Case 1.}
Assume for a nontrivial critical configuration,
$\bar x_1=0$, $\bar x_2\ne0$, and $\bar x_3\ne0$.
Define $I_1=\{ i: \alpha_i=\lambda_1\}$. We claim $|I_1|=2$.

Let $I_0$ be the complement of $I_1$.
First observe that for $i\in I_0$ we have $x_{i1}=0$ because
$(\alpha_i-\lambda_1)x_{i1}=\bar x_1=0$.
Then $I_1\ne\emptyset$, otherwise the configuration
lies in the plane $x_1=0$, hence trivial.
Having exactly one nonzero component among $x_{i1}$
is not possible either because $\bar x_1=0$.
Hence $|I_1|\ne1$.

Consider $|I_1|\ge 3$.
Since $\bar x_k\ne0$ for $k=2,3$, then
$x_{ik}=\frac{\bar x_k}{\alpha_i-\lambda_k}$ for $k=2,3$.
Then $x_{ik}$ is independent of $i\in I_1$ for $k=2,3$,
so the points $p_i$ for $i\in I_1$ can differ only in the first
component $x_{i1}$. But since $|p_i|=1$, then the
only remaining freedom is in the sign of $x_{i1}$.
So if $|I_1|\ge 3$, then the points have to repeat.

Hence $|I_1|=2$, the two points $p_i$, $i\in I_1$,
are symmetric about the plane $x_1=0$,
and the three points $p_i$, $i\in I_0$,
are in the plane $x_1=0$.
$\square$
\medskip

{\bf Case 2.}
Let $\bar x_1=0$, $\bar x_2=0$, and $\bar x_3\ne0$.
Define $I_1$ and $I_0$ as above.

Subcase $\lambda_1\ne\lambda_2$.
If $i\in I_0$, then $x_{i1}=0$.
If $i\in I_1$, then $x_{i2}=0$ because $\alpha_i\ne\lambda_2$.
Then $|I_1|=2$ by the same reason as in Case 1,
and the configuration falls into Main Special Case.

Subcase $\lambda_1=\lambda_2=:\lambda$.
If $i\in I_0$, then $x_{i1}=x_{i2}=0$, and $x_{i3}=\pm 1$.
If $i\in I_1$, then $x_{i3}=\frac{\bar x_3}{\lambda-\lambda_3}$
is independent of $i$.

If $I_0=\emptyset$, then all the points $p_i$ lie
on the same circle, which we regard as trivial.

If $|I_0|=1$, then one point is, say $(0,0,1)$,
the North pole, and the rest are on the same parallel.
Such a configuration will have to be \FP, a regular pyramid
with square base, which falls into Main Special Case because
it is symmetric about the plane through the vertex and
a diagonal of the base.

If $|I_0|=2$, then two points are $(0,0,\pm1)$,
and the rest lie on the same non-equatorial parallel.
Such a configuration is not critical.

The case $|I_0|>2$ can not occur because
otherwise the points will repeat.
$\square$
\medskip

{\bf Case 3.}
We have $\bar x_1=\bar x_2=\bar x_3=0$.
Define
$I_k=\{ i: \alpha_i=\lambda_k\}$, $k=1,2,3$.
Note that every index $1\le i\le 5$ has to be in one of $I_k$,
otherwise $p_i=0$.

Subcase $\lambda_k$, $1\le k\le 3$, are distinct.
By the same arguments as above we obtain that
either $I_1=\emptyset$ or $|I_1|=2$ and
the corresponding points $p_i$ are $(\pm1,0,0)$.
Similarly, for all $k$ the cardinalities
$|I_k|$ can be only 0 or 2, which is not possible because
5 is an odd number.

Subcase $\lambda_1=\lambda_2\ne\lambda_3$.
Then for $i\in I_1=I_2$, we have $x_{i3}=0$.
Also $|I_3|=2$, and for $i\in I_3$,
we have $p_i=(0,0,\pm1)$.
This configuration will have to be \TBP, and it falls
into Main Special Case.

Finally, the subcase $\lambda_1=\lambda_2=\lambda_3$
can not occur, which we prove below.
$\square$
\medskip

\begin{prop}
\label{equal-lambdas}
For $m=3$, $n=5$ there is no configuration
for which $\bar x=0$ and $\lambda_1=\lambda_2=\lambda_3$.
\end{prop}

\proof
Suppose such a configuration does exist.
Let $X$ be the matrix of the configuration.
Then the matrix $Y=X/\sqrt{\lambda_1}$ has orthonormal
columns. Add a fourth column to $Y$ with entries all
equal to $1/\sqrt{5}$. Since $\bar x=0$, then the resulting
matrix still has orthonormal columns. Add the fifth
column $(y_{i5})$ to obtain an orthogonal matrix.
Note the rows of an orthogonal matrix are unit vectors.
Since all $|p_i|=1$, and the lambdas are equal,
then $|y_{i5}|$ will have to be
equal for all $i$, in fact $y_{i5}=\pm1/\sqrt{5}$.
This column also must be orthogonal to the fourth column,
which has equal entries, hence $\sum y_{i5}=0$.
However, the five terms $\pm1/\sqrt{5}$ can not add up
to zero.
$\square$
\medskip

\begin{remark}{\rm
It turns out that for $m=3, n>5$ the case $\bar x=0$,
$\lambda_1=\lambda_2=\lambda_3$ can occur and provides
minimizers for the biquadratic energy.
The minimizers form continuous families.
}\end{remark}

\section{Main Special Case}

While the $\alpha$-$\lambda$-notation proved to be
useful in obtaining geometric restrictions on
the critical configurations, we were unable to use it
effectively in solving the equilibrium equations
in this case. We do it here by brute force.

A configuration $P$ in Main Special Case has the form
\begin{align*}
&p_i=(\cos\t_i,\sin\t_i,0),  &&1\le i\le 3, \\
&p_\alpha=(r,0,\pm\sqrt{1-r^2}),  &&\alpha=4,5.
\end{align*}
Here $\t_i$, $1\le i\le 3$, and $|r|<1$
are real parameters.
For $r=0$ it is immediate that the only nontrivial
critical configuration is $\TBP$.
Hence we can assume $r\ne 0$. We have
\begin{align*}
&b_{ij}=p_i\cdot p_j=\cos\t_{ij}, &\t_{ij}:=\t_i-\t_j,\\
&b_{i\alpha}=r\cos\t_i,\\
&b_{45}=2r^2-1.
\end{align*}
Instead of using the general equilibrium equations
\eqref{Crit4},
we derive them anew in this case.
The biquadratic energy has the form
\[
E=E(P)=\sum_{i<j\le5}(b_{ij}+1)^2
=\sum_{i<j\le3}(\cos\t_{ij}+1)^2
+2\sum_{i=1}^3(r\cos\t_i+1)^2+4r^4.
\]
Then the equations $\partial E/\partial\t_i=0$ and
$\partial E/\partial r=0$ yield
\begin{align*}
&\sum_{j=1}^3(\cos\t_{ij}+1)\sin\t_{ij}
+2(r\cos\t_i+1)r\sin\t_i=0,  &1\le i\le 3, \\
&\sum_{j=1}^3(r\cos\t_j+1)\cos\t_j+4r^3=0.
\end{align*}
We rewrite the above equations in terms of the complex variables
$z_j=e^{i\t_j}$.
\begin{align}
&\sum_{j=1}^3\left(\frac{z_i^2}{z_j^2}-\frac{z_j^2}{z_i^2}\right)
+2\sum_{j=1}^3\left(\frac{z_i}{z_j}-\frac{z_j}{z_i}\right)
+2r^2(z_i^2-z_i^{-2})
+4r(z_i-z_i^{-1})=0,  &1\le i\le 3, \label{eq1}\\
&r\sum_{j=1}^3(z_j^2+z_j^{-2})
+2\sum_{j=1}^3(z_j+z_j^{-1})
+16r^3+6r=0. \label{eq2}
\end{align}
Adding up the three equations \eqref{eq1},
assuming $r\ne0$, we obtain
\[
r\sum(z_j^2-z_j^{-2})
+2\sum(z_j-z_j^{-1})=0.
\]
The latter and \eqref{eq2} yield
\begin{equation}
\label{eq2-1}
r\sum z_j^2 +2\sum z_j
=r\sum z_j^{-2} +2\sum z_j^{-1}
=-8r^3-3r.
\end{equation}
Introduce the power sums
$s_p=\sum z_j^p$, $p=\pm1, \pm2$.
We {\it formally} rewrite each equation in \eqref{eq1}
as a linear combination of powers of $z_i$ with
coefficients depending on $z_j$-s, that is,
\begin{equation}
\label{eq3}
Az_i^2+Bz_i+Cz_i^{-1}+Dz_i^{-2}=0.
\end{equation}
Then the coefficients have the form
\begin{equation}
\label{eq4}
A=s_{-2}+2r^2,\quad
B=2s_{-1}+4r,\quad
C=-2s_1-4r,\quad
D=-s_2-2r^2.
\end{equation}
The equation \eqref{eq3} is non-trivial:
the case $A=B=C=D=0$ can not occur at all,
and if $A=0$, $|z_i|=1$, then the numbers $z_i$
can not be distinct -
we leave the details to the reader.

Let $\sigma_p$, $1\le p\le3$, be elementary symmetric
functions of $z_i$.
Then the three distinct numbers $z_i$ satisfy both
\eqref{eq3} and
$z^3-\sigma_1 z^2+\sigma_2 z-\sigma_3=0$.
Let $z_0$ be the fourth root of \eqref{eq3},
which in principle may coincide with one of $z_i$.
Then we have the identity
\begin{equation}
\label{eq5}
Az^4+Bz^3+Cz+D
=A(z-z_0)(z^3-\sigma_1 z^2+\sigma_2 z-\sigma_3).
\end{equation}
The goal of these manipulations is converting the system
(\ref{eq1}--\ref{eq2}) to the new variables
$x=s_1$, $y=s_{-1}$, $t=\sigma_3$.
We have
\begin{equation}
\label{eq6}
s_1=x,\quad
s_2=x^2-2yt,\quad
s_{-1}=y,\quad
s_{-2}=y^2-2x/t,\quad
\sigma_1=x,\quad
\sigma_2=yt,\quad
\sigma_3=t.
\end{equation}
Equating the coefficients of powers of $z$ in \eqref{eq5}
yields
\begin{align}
&B=2y+4r=-A(z_0+x),
&&0=A(z_0x+yt),\label{eq7}\\
&C=-2x-4r=-A(z_0yt+t),
&&D=-x^2+2yt-2r^2=Atz_0.\label{eq8}
\end{align}
Eliminating $z_0=-yt/x$ yields
\begin{align}
A(-yt+x^2)+2xy+4rx=0,\label{eq10}\\
A(-y^2t^2+xt)-2x^2-4rx=0,\label{eq11}\\
Ayt^2-x^3+2xyt-2r^2x=0.\label{eq12}
\end{align}
We modify \eqref{eq10} by multiplying by $t$
and adding \eqref{eq12} to it.
We multiply \eqref{eq12} by $y$ and add to
\eqref{eq11}. We plug $A=y^2-2x/t+2r^2$
in \eqref{eq12}. Finally, we solve each equation
for $t$ and obtain
\begin{equation}
\label{eq-t}
t=\frac{3x^2+2r^2}{xy^2+2r^2x+4y+4r}
=\frac{x^2y+2r^2y+4x+4r}{3y^2+2r^2},\quad
t^2=\frac{x(x^2+2r^2)}{y(y^2+2r^2)}.
\end{equation}
We note that neither numerators nor denominators
in \eqref{eq-t} can vanish.
We already discarded $r=0$.
The case $x=0$ can not occur because then $y=\bar x=0$,
and by (\ref{eq7}-\ref{eq8}) we eventually get $r=\pm 2$,
which is not allowed. The other numerators or denominators
can not vanish either by similar reasons - we leave the details
to the reader.

Eliminating $t$ yields
\begin{align}
&(3x^2+2r^2)(3y^2+2r^2)-(x^2y+2r^2y+4x+4r)(xy^2+2r^2x+4y+4r)=0,
\label{eq13}\\
&x(x^2+2r^2)(3y^2+2r^2)(xy^2+2r^2x+4y+4r)\notag\\
&\qquad\qquad\qquad\qquad\qquad
-y(y^2+2r^2)(3x^2+2r^2)(x^2y+2r^2y+4x+4r)=0.
\label{eq14}
\end{align}
The equations (\ref{eq13}--\ref{eq14})
reflect only \eqref{eq1}.
Converting \eqref{eq2-1} to the new variables yields
\begin{equation}
(3x^2+2r^2)(ry^2+2y+8r^3+3r)-
2rx(xy^2+2r^2x+4y+4r)=0,
\label{eq15}
\end{equation}
and the one with the interchange $x\leftrightarrow y$.
The system (\ref{eq13}--\ref{eq15}) is simple enough
for symbolic computations software.
We consider two cases $x\ne y$ and $x=y$.

Let $x\ne y$. Note \eqref{eq14} is antisymmetric in
$x$ and $y$. Then dividing \eqref{eq14} by $r(x-y)$ yields
\begin{align}
8 r^4 + 4 r^5 x + 4 r^2 x^2 + 2 r^3 x^3 + 4 r^5 y - 8 r^2 x y -
 8 r x^2 y + 2 r^3 x^2 y + 4 r^2 y^2 \notag\\
- 8 r x y^2 + 2 r^3 x y^2 +
 6 x^2 y^2 + r x^3 y^2 + 2 r^3 y^3 + r x^2 y^3=0.
\label{eq16}
\end{align}
Eliminating $x$ and $y$ from (\ref{eq13}, \ref{eq15}, \ref{eq16})
and (\ref{eq15}) with the interchange $x\leftrightarrow y$
using Macaulay2 yields
\begin{align*}
&7744275\, r^4 + 80139015\, r^6 + 351783930\, r^8 + 861239064\, r^{10}
+ 1282072196\, r^{12} + 1176047932\, r^{14} \\
&+ 632113944\, r^{16} + 172153584\, r^{18} + 18541440\, r^{20}
+ 3882816\, r^{22} + 777600\, r^{24} + 186624\, r^{26}=0.
\end{align*}
Since all the coefficients are positive, then the only real
solution is the trivial $r=0$.

We now consider the case $x=y$. Recall $y=\bar x$.
Then $x=y$ implies $x\in\R$.
The last equation in \eqref{eq-t} implies $t=\pm1$.
Note that the system in either $(z_i,r)$ or $(x,y,t,r)$
is invariant under changing the sign of all the variables,
so without loss of generality $t=1$.
We have $\sigma_1=\sigma_2=x$, $\sigma_3=1$. Then $z_i$-s
are solutions of the equation
$z^3-xz^2+xz-1=(z-1)(z^2+(1-x)z+1)=0$.
Then, say $z_1=1$, $\Re z_2=\Re z_3=(x-1)/2$.

For $x=y$ the equation \eqref{eq14} holds automatically,
and the equations (\ref{eq13}, \ref{eq15})
for $t=1$ reduce to
\begin{align}
&(3x^2+2r^2)-(x^3+2r^2x+4x+4r)=0,
\label{eq20}\\
&(rx^2+2x+8r^3+3r)-2rx=0.
\label{eq21}
\end{align}
Eliminating $x$ yields
\begin{equation*}
r(1 + 2 r)(1 + 2 r + 2 r^2)(8 - 9 r + 6 r^2)(1 + 3 r + 6 r^3)=0.
\end{equation*}
The real nonzero solutions are $r=-1/2$ and the
solutions of the equation
\begin{equation}
\label{eq22}
6r^3+3r+1=0.
\end{equation}

Let $r=-1/2$. Then one can find $x=1$, $\Re z_2=\Re z_3=0$,
$z_2,z_3=\pm i$. This configuration is $\TBP$.

The equation \eqref{eq22} has only one real root
$r\approx-0.286$. Eliminating $r$ from
(\ref{eq20}-\ref{eq22}) yields $3x^3-9x^2+15x-5=0$.
The latter matches with \eqref{eq22} if $r=(x-1)/2$.
Then the corresponding configuration consists of
$(1,0,0)$, $(r, \pm\sqrt{1-r^2}, 0)$, and
$(r,0,\pm\sqrt{1-r^2})$. This is a four-pyramid with
square base \FP. The inner products have the following values:
four $r$, four $r^2$, and two $2r^2-1$.
Then the energy is equal to
$E(\FP)=4(r+1)^2+4(r^2+1)^2+2(2r^2)^2\approx 7.9$,
which is bigger than $E(\TBP)=6.75$.
Hence $\TBP$ is a unique minimizer.
The proof of Theorem \ref{MainThm} is complete.


{\footnotesize

\begin{thebibliography}{CIT}

\bibitem{CK}
H. Cohn and A. Kumar,
{\it Universally optimal distributions of points on spheres},
J. Amer. Math. Soc. 20 (2007), 99--148.

\bibitem{DLT}
P. D. Dragnev, D. A. Legg, and D. W. Townsend,
{\it Discrete logarithmic energy on the sphere},
Pacific J. Math. 207 (2002), 345--357.

\bibitem{HouShao}
X. Hou and J. Shao,
{\it Spherical distribution of 5 points with
maximal distance sum},
arXiv: 0906.0937, 45 pp.

\bibitem{Sch} R. E. Schwartz,
{\it The 5 electron case of Thomson's problem},
arXiv: 1001.3702, 67 pp.

\end{thebibliography}
}
\end{document}